\documentclass[review]{elsarticle}

\usepackage{hyperref}

\journal{Ecological Modelling}

\usepackage{numcompress}\biboptions{authoryear}

\usepackage{amsmath}
\usepackage{amssymb}
\usepackage{amsthm}

\newdefinition{remark}[theorem]{Remark}

\newdefinition{genvor}{General assumption}[section]

\newdefinition{definition}[theorem]{Definition}
\newdefinition{erg}[theorem]{Result}

\renewcommand{\vec}[1]{\mathbf{#1}}

\begin{document}

\begin{frontmatter}

\title{Mathematical analysis of the $PO_4$-$DOP$-$Fe$ marine ecosystem model driven by 3-D ocean transport
}


\author[mymainaddress]{Christina Roschat\corref{mycorrespondingauthor}}
\cortext[mycorrespondingauthor]{Corresponding author. Tel.: +49 431 880 7452}
\ead{cro@informatik.uni-kiel.de}

\author[mymainaddress]{Thomas Slawig}

\address[mymainaddress]{Dep. of Computer Science, Kiel Marine Science - Centre for Interdisciplinary Marine Science, Christian-Albrechts-Universit\"at zu Kiel, 24098 Kiel}

\begin{abstract}
Marine ecosystem models are developed to understand and simulate the biogeochemical processes involved in marine ecosystems. Parekh, Follows and Boyle introduced the $PO_4$-$DOP$-$Fe$ model of the coupled phosphorus and iron cycles in 2005. Especially the part describing the phosphorus cycle ($PO_4$-$DOP$ model) is often applied in the context of 
parameter identification.
The mathematical analysis presented in this study is concerned with the existence of solutions and the reconstruction of parameters from given data. Both are important questions in the numerical model's assessment and validation not answered so far.
%
In this study, we obtain transient, stationary and periodic solutions (steady annual cycles) of the $PO_4$-$DOP$-$Fe$ model equations after a slight change in the equation modeling iron. 
This result confirms the validity of the solutions computed numerically. 
Furthermore, we present a calculation showing that four of the $PO_4$-$DOP$ model's parameters are possibly dependent, i.e. different parameter values might be associated with the same model output. Thereby, we identify a relevant source of uncertainty in parameter identification. On the basis of the results, possible ways to overcome this deficit can be proposed. 
In addition, the stated mathematical conditions for solvability are universal 
and thus applicable to the analysis of other ecosystem models as well.
 \end{abstract}


\end{frontmatter}


\section{Introduction}
Being a part of the global carbon cycle marine ecosystems considerably influence on the earth's climate. 
Marine ecosystem models, describing the biogeochemical processes involved, provide an important tool to understand these processes and thus to predict the future concentration of carbon dioxide in oceans and atmosphere.

One important example for models of this kind is the $PO_4$-$DOP$-$Fe$ model by \cite{pa05}, describing the concentrations of three tracers phosphate ($PO_4$), dissolved organic phosphorus ($DOP$) and iron ($Fe$). Following \cite{kri10}, we use the synonymous name $N$-$DOP$-$Fe$ model. The abbreviation $N$ stands for the more general expression ``nutrient''.

Fixing the iron concentration, the evolving $N$-$DOP$ model, consisting of the first two tracers, describes the marine phosphorus cycle. Despite the low complexity and associated simplifications, it is still in use \citep{kri12,pa06}. 
An assessment on the basis of oceanic observations even indicates that, under certain circumstances, the $N$-$DOP$ model can compete with more complicated ones \citep{kri10, kri12}. 
Due to the low complexity, it is often used for the purpose of testing numerical methods, e.g. in \citet{pr13}. 

The third tracer expresses the influence of iron on the phosphorus cycle. Iron is an important nutrient. For instance, \cite{kri12} assume that their observed misfit of model output and observational data might be due to missing iron limitation. The model can also picture the changes in the phosphorus cycle induced by artificial fertilization with iron.

The $N$-$DOP$-$Fe$ model consists of three advection-diffusion-reaction equations, each characterizing one tracer concentration on a three-dimensional ocean domain. The concentrations are influenced by ocean transport, i.e. advection and diffusion, and biogeochemical processes modeled by specific reaction terms.

The first part of this work studies the existence of transient, periodic and stationary solutions of the $N$-$DOP$-$Fe$ model equations. Periodic solutions (steady annual cycles) are the most relevant ones. Such results are desirable because the numerical model output used in applications is computed by approximating a solution of the original model equations. Therefore, the quality and reliability of the model output depends on the solvability of these equations. A comparable analysis of the $N$-$DOP$-$Fe$ model equations has not been undertaken so far.
%
%

Any model has to be calibrated, i.e. adapted to the ecosystem in question. This is achieved by the choice of the model's parameters. Parameters, like e.g. remineralization rates or half saturation constants, are essential quantities characterizing the processes modeled. The $N$-$DOP$ model includes seven parameters. 

Parameters can be estimated by means of laboratory experiments. If such measurements are too difficult or expensive, parameters are alternatively identified via optimization with respect to observational data. Here, a parameter set minimizing the difference between model output and data is determined. An important area of application of the $N$-$DOP$ model is the testing of numerical methods designed to solve such kind of minimization problems \citep{pr13}.
The method in question is applied to synthetic data, i.e. data corresponding to known optimal parameters, assuming that a correct method is able to identify these.
However, this assumption only holds true if all parameters are uniquely identifiable, i.e. if each possible model output is associated with a single parameter set. 

The second part of this work is dedicated to the question, unanswered so far, which of the $N$-$DOP$ model parameters are uniquely identifiable. We name and justify different alternatives to alter the reaction terms in such a way that all parameters become identifiable.

The mathe\-matical results about the $N$-$DOP$-$Fe$ model are introduced and explained to readers from other disciplines than mathematics. By outlining the mathematical proceeding, we intend to provide an overview about which mathematical conditions and assumptions are responsible for the model's properties. These may be a guideline for the investigation and development of other models or alternative reaction terms. 

The paper is structured as follows: In the next section, we introduce the classical and the weak formulation of the $N$-$DOP$-$Fe$ model equations, the latter being the object of investigation. The next three sections each deal with one type of solution. The corresponding results concerning existence and, where possible, uniqueness are formulated and shortly justified. In Sec.~\ref{sec:paropt}, we investigate identifiability of parameters. In the last two sections, the results are discussed and conclusions are drawn.

\section{Model equations}
The following introduction of the $N$-$DOP$-$Fe$ model is based on the original work by~\citet{pa05} 
as well as on \citet{ro14,ro214}. 

The ecosystem is located in a three-dimensional bounded domain $\Omega$ determined by the bounded water surface $\Omega^{\prime}\subseteq\mathbb{R}^2$ and the depth $h(x^{\prime})>0$ at every surface point $x^{\prime}\in \Omega^{\prime}$. 
The boundary $\Gamma$ is the union of the surface $\Gamma^{\prime}$ and the boundary inside the water
. 

The domain is separated into two layers, the euphotic, light-flooded zone $\Omega_1$ up to a depth of $\bar{h}_e:=120\,\text{m}$ and the dark, aphotic zone $\Omega_2$ beneath. 
The actual depth of the euphotic zone beneath some surface point $x^{\prime}$ is defined by $h_e(x^{\prime}):=\min\{\bar{h}_e,h(x^{\prime})\}$. We split the surface into the part $\Omega_2^{\prime}:=\{x^{\prime}\in \Omega^{\prime}; h(x^{\prime})>\bar{h}_e\}$ above the aphotic zone and the rest $\Omega^{\prime}_1:=\Omega^{\prime}\setminus \Omega_2^{\prime}$. The boundary is analogously divided into the euphotic part $\Gamma_1$ and the aphotic part $\Gamma_2$.
%
%

For the sake of a clearer terminology, we write henceforth $y_1$ for the phosphate concentration, $y_2$ for the concentration of $DOP$ and $y_3$ for the iron concentration. The three tracers, assembled in the vector $\vec{y}=(y_1,y_2,y_3)$, are characterized by the system of advection-diffusion-reaction equations 
 \begin{equation}
\label{eq:zustand}
\begin{array}{rcll}
\displaystyle \partial_ty_1+\textrm{div}(\vec{v}y_1)-\textrm{div}(\kappa\nabla  y_1)-\lambda y_2+d_1(\vec{y})&=&0
\\
\displaystyle \partial_ty_2+\textrm{div}(\vec{v}y_2)-\textrm{div}(\kappa\nabla  y_2)+\lambda y_2+d_2(\vec{y})&=&0
\\
\displaystyle \partial_ty_3+\textrm{div}(\vec{v}y_3)-\textrm{div}(\kappa\nabla  y_3)+J_{Fe}-\lambda y_2R_{Fe}+d_3(\vec{y})&=&S_{Fe},
\end{array}
\end{equation}
each dependent on space and time in a finite time span $[0,T]$.  Regarding periodic solutions, a reasonable value for the final point of time $T$ is one year.

%
%
The second and third terms on the left-hand sides describe the ocean dynamics by means of the velocity field $\vec{v}$ and the diffusion coefficient $\kappa$. Sometimes these terms are summarized in a linear operator or, if already discretized, in a matrix \citep{kri10, kri12}. 
To reduce computational effort, the models are often run in an ``offline'' mode, i.e. the influence the tracers have on the oceans dynamics is neglected. This has the advantage that the values of $\vec{v}$ and $\kappa$ can be precomputed with one run of an ocean circulation model. For that reason, $\vec{v}$ and $\kappa$ 
are assumed to be known in the theoretical analysis as well.  
Since turbulent exceeds molecular diffusion by far the values of $\kappa$ are assumed equal in all equations.

Furthermore, since the ocean is a closed system, it is reasonable to assume that the velocity field $\vec{v}$ is divergence free and does not point out of the boundaries. These two properties are crucial for the mathematical analysis. 

The reaction terms $d_j(\vec{y}),J_{Fe},\lambda y_2$ and $\lambda y_2R_{Fe}$ describe the biogeochemical coupling. Therefore, they depend on one or several of the tracers. On the right-hand side of each equation, sources and sinks of the respective tracer are displayed. Only iron has a non-zero source term. 

In the following two subsections, we introduce and derive the reaction terms associated with the phosphorus cycle and the iron equation, respectively.
Two further subsections deal with boundary conditions and weak solutions.

\subsection{The phosphorus cycle}\label{sec:phosphor}
%
%
One important process is the remineralization of $y_2$ into $y_1$. In the first two equations of Eq.~\eqref{eq:zustand}, it is modeled as a first order loss process with a remineralization rate $\lambda$ between 0 and 1. Being independent of light, this transformation takes place in the whole domain $\Omega$. The remaining processes, represented by the reaction terms $d_1$ and $d_2$, differ according to the layers. In the euphotic zone, $y_1$ is consumed by the photosynthesis of phytoplankton. This uptake is modeled 
by
\begin{align*}
G(y_1,y_3)&:=
\alpha\frac{y_1}{|y_1|+K_P}\frac{y_3}{|y_3|+K_F}\frac{I\text{e}^{-x_3K_W}}{|I\text{e}^{-x_3K_W}|+K_I}. 
\end{align*}
Here, the maximum uptake $\alpha>0$ is limited by the present concentrations of phosphate $y_1$ and iron $y_3$ as well as insolation via Michaelis-Menten kinetics. The corresponding saturation functions are equipped with the half saturation constants $K_P,K_F,K_I$. The absolute values in the denominators ensure that the fractions are mathematically well-defined because it is a priori unknown if the solutions $y_1$ and $y_3$ are nonnegative. If this holds true, the above formulation is in accordance with Michaelis-Menten kinetics.
Incidence of light is formalized by the insolation $I=I(x^{\prime},t)$ depending on the water surface and time. The incident light decreases exponentially with depth $x_3$ and the attenuation coefficient $K_W$ for seawater. Below the euphotic zone, the values of $I$ are zero \citep{palpla76}. 

A fraction $\nu\in(0,1]$ of the uptake $G$ is transformed into $y_2$ while the remnants, integrated over the whole water column, are exported into the aphotic zone $\Omega_2$. The remineralization during the sinking of particles is described by a parameter $b$. 

The outlined processes are represented by nonlinear coupling terms with different appearance in each layer. In the euphotic zone $\Omega_1$, they are given by 
\begin{align*}
d_1(\vec{y})&:= G(y_1,y_3)\\
d_2(\vec{y})&:= -\nu G(y_1,y_3),
\end{align*}
and, in $\Omega_2$, by
\begin{align*}
d_1(\vec{y})&:=-(1-\nu)\int_0^{h_e}G(y_1,y_3)dx_3
\frac{b}{\bar{h}_e}\left(\frac{x_3}{\bar{h}_e}\right)^{-b-1}\\
d_2(\vec{y})&:=0.
\end{align*}
We remark that the $N$-$DOP$ model contains the seven parameters $\lambda,\alpha,K_{P},\\K_{I},K_{W},b,\nu$.

%
%

\subsection{The iron equation}\label{subsec:iron}
In opposite to phosphorus, there is a source term for iron given by $S_{Fe}:=\beta F_{in}$ which is non-zero only in the euphotic zone. The parameter $\beta$ represents the solubility of iron in seawater and $F_{in}$ quantifies the aeolian source of iron. 

The reaction terms express how the phosphorus cycle, complexation and scavenging influence on the iron cycle. The iron concentration increases with re\-mineralization and decreases with consumption of phosphorus. Being expressed in phosphorus units, these values are multiplied by the constant ratio $R_{Fe}
$ in order to become iron units. Thus, iron increases with $\lambda y_2R_{Fe}$. The decrease induced by consumption is expressed by the coupling term 
\begin{displaymath}
d_3(\vec{y})=G(y_1,y_3)R_{Fe}.
\end{displaymath}

\subsubsection{Scavenging and complexation - original formulation}
The summand $J_{Fe}$ represents the influence of complexation and scavenging on the iron concentration. Since a certain amount of iron is complexed with organic ligand, the total iron $y_3$ is split into free iron $Fe^{\prime}$ and complexed iron $FeL$. Similarly, the total ligand $L_T$ is the sum of free ligand $L^{\prime}$ and the complexed ligand, equal to the complexed iron $FeL$. These relations are expressed by the  formulae $y_3=Fe^{\prime}+FeL$ and $L_T=L^{\prime}+FeL$.
\cite{pa05} set $L_T=1$. 

Only the free iron is subject to scavenging which is thus modeled as the first order loss process $J_{Fe}:=\tau k_0C_p^{\Phi}Fe^{\prime}$. Here, $\tau,\Phi$ are numbers, $k_0$ is the initial scavenging rate and $C_p$ represents the particle concentration which decreases with depth. 

Since a feasible mathematical formulation requires that $J_{Fe}$ depends at least on one of the tracers $y_1,y_2$ or $y_3$, we express $Fe^{\prime}$ using $y_3$. To this end, \citet{pa05} additionally provide the equilibrium relationship $K=FeL/Fe^{\prime}L^{\prime}$ with a positive constant $K$. Inserting the equivalent expression $L^{\prime}=FeL/KFe^{\prime}$ into the equation for ligand we obtain
\begin{displaymath}
L_T=FeL+L^{\prime}=FeL+\frac{FeL}{KFe^{\prime}}=FeL\left(1+\frac{1}{KFe^{\prime}}\right).
\end{displaymath}
$FeL$ can be replaced by $y_3-Fe^{\prime}$. This gives
\begin{displaymath}
L_T=(y_3-Fe^{\prime})\left(1+\frac{1}{KFe^{\prime}}\right)=y_3-Fe^{\prime}+\frac{y_3}{KFe^{\prime}}-\frac{1}{K}.
\end{displaymath}
With the abbreviation $H(y_3):=L_T+1/K-y_3$, this proves equivalent to 
\begin{displaymath}
Fe^{\prime\,2}+H(y_3)Fe^{\prime}-\frac{y_3}{K}=0.
\end{displaymath}
The two solutions of this quadratic equation are known. It can be easily shown that one of them has only negative values and is therefore inappropriate to describe the amount of free iron. Thus, we find
\begin{displaymath}
Fe^{\prime}(y_3)=-\frac{1}{2}H(y_3)+\sqrt{\frac{H(y_3)^2}{4}+\frac{y_3}{K}}.
\end{displaymath}
To ensure that the square root is real we show that the radicand $r:=(L_T+1/K-y_3)^2/4+{y_3}/{K}$ is positive. This is obvious whenever $y_3$ is positive, since in this case $r\geq{y_3}/{K}>0$. In the non-positive case, 
we transform $r$ into
\begin{align*}
r &=\frac{1}{4}\left(L_T+\frac{1}{K}-y_3\right)^2+\frac{y_3}{K}\\
 &=
 \frac{1}{4}\left(L_T+\frac{1}{K}\right)^2-\frac{1}{2}\left(L_T+\frac{1}{K}\right)y_3+\frac{1}{4}y_3^2+\frac{y_3}{K}\\
&=
 \frac{1}{4}\left(L_T+\frac{1}{K}\right)^2+\frac{1}{4}y_3^2-\frac{1}{2}y_3\left(L_T-\frac{1}{K}\right).
\end{align*}
Since $K$ tends to be a large number and $L_T$ takes a value around one (\cite{pa05} set $K=\exp(11)$ and $L_T=1$) the expression $L_T-1/K$ is nonnegative. Therefore, $y_3\leq0$ yields $-y_3(L_T-1/{K})/2\geq0$ and, as a consequence, $r\geq(L_T+{1}/{K})^2/4>0$. In particular, this result implicates that $Fe^{\prime}$ is 
differentiable everywhere on $\mathbb{R}$. It can be shown that the derivative is positive and thus that $Fe^{\prime}$ is a monotonically increasing real function. 

\begin{figure}[t]
\includegraphics[width=8.3cm]{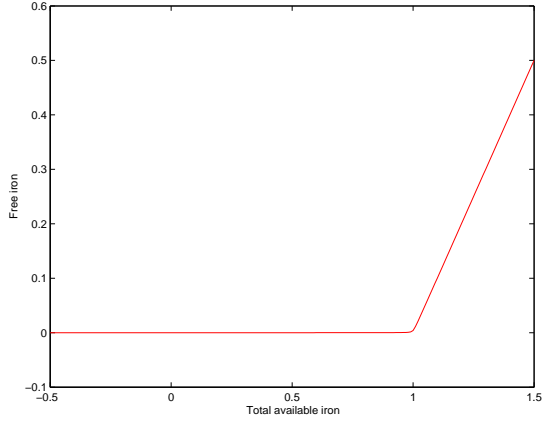}
\caption{Graph of free iron $Fe^{\prime}$ in dependence of total available iron $y_3$ assuming $K=\exp(11)$ and $L_T=1$ \citep{pa05}.}
\label{fig:funkJ}
\end{figure}
In Fig.~\ref{fig:funkJ}, we see that, as long as approximately $y_3<L_T=1$, the values for $Fe^{\prime}$ increase very slowly, for $y_1>L_T$ they increase with a gradient of almost one. This behavior corresponds to the statement of \cite{pa05} that they ``rapidly precipitate $Fe^{\prime}$ when $Fe_T>L_T$'' meaning that, as long as ligand is available, only a small amount of iron remains free.

Thus, the reaction term associated with scavenging and complexation depending on the third tracer $y_3$ is determined by
\begin{displaymath}
J_{Fe}(y_3):=\tau k_0C_p^{\Phi}Fe^{\prime}(y_3).
\end{displaymath}

%
\subsubsection{Scavenging and complexation - adjusted formulation}
The reaction term $J_{Fe}(y_3)$ derived above lacks one important property with regard to solvability. This is due to the fact that the function 
\begin{displaymath}
Fe^{\prime}(y_3)=-\frac{1}{2}\left(L_T+\frac{1}{K}-y_3\right)+\sqrt{\frac{1}{4}\left(L_T+\frac{1}{K}-y_3\right)+\frac{y_3}{K}}
\end{displaymath}
modeling free iron levels off fast in negative direction.
For that reason, we propose an alternative term $FeF(y_3)$ for free iron. Fig.~\ref{fig:funkJ} indicates that $Fe^{\prime}$ resembles a straight line with a slope of 1 for approximately $y_3>L_T$.
 If $y_3\leq L_T$, the curve tends to zero very fast. As a substitute for $Fe^{\prime}$ we therefore define a piecewise linear function, composed of a line with a slope of 1 and another line through zero. As their ``meeting point'' we determine $(L_T,Fe^{\prime}(L_T))$ where the function $Fe^{\prime}$ apparently turns upwards.
Formally, the function $FeF$ depending on total iron $y_3$ is defined by
\begin{displaymath}
FeF(y_3)=
\begin{cases}
 y_3+Fe^{\prime}(L_T)-L_T & \text{if $y_3>L_T$}, \\
\frac{Fe^{\prime}(L_T)}{L_T}y_3   & \text{if $y_3\leq L_T$}.
\end{cases}
\end{displaymath}
\begin{figure}
\includegraphics[width=8.3cm]{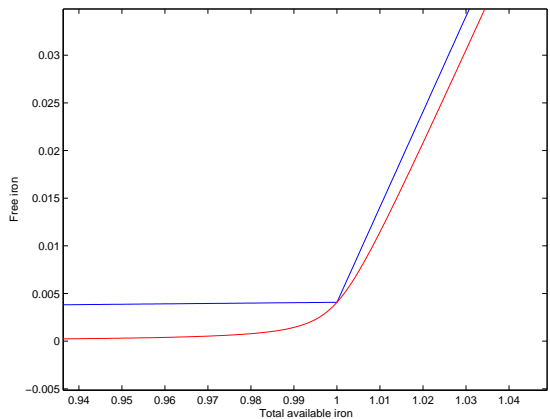}
\caption{Graph of free iron $Fe^{\prime}$ (red line) in comparison to the alternative function $FeF$ (blue line) around $y_3=L_T$ assuming $K=\exp(11)$ and $L_T=1$ \citep{pa05}.}
\label{fig:funklin}
\end{figure}

Figure~\ref{fig:funklin} shows the difference between exemplary curves of $Fe^{\prime}$ and $FeF$ around $y_3=L_T$. Comparing the complete curves reveals that, in the positive region, the piecewise linear function $FeF$ lies slightly above the original and their distance remains small. In the negative region, the curves behave conversely. The distance increases with decreasing values of $y_3$ because $FeF$ has a constant slope and $Fe^{\prime}$ a decreasing one. However, negative values for $y_3$ are not relevant in most applications. 

Taking into account these considerations, a possible alternative for $J_{Fe}$ might be the adjusted reaction term
\begin{displaymath}
 J(y_3):=\tau k_0C_p^{\Phi}FeF(y_3).
\end{displaymath}

%
\subsection{Boundary conditions} 
The original $N$-$DOP$-$Fe$ model formulation lacks explicit statements about the tracers' behavior on the boundary. As this information is needed 
with respect to a weak formulation of the model equations (cf. Sec.~\ref{sec:weak}), suitable conditions are derived in this section. 
In general, a boundary condition of Neumann type for the $j$-th tracer has the form
\begin{equation}
\label{eq:boundarycondition}
\nabla y_j\cdot(\kappa\eta)+ b_j(\vec{y})=0.
\end{equation} 
The symbol $\eta$ stands for the outward-pointing unit normal vector field on the boundary. Therefore, $\nabla y_j\cdot(\kappa\eta)$ can be understood as the change of concentration alongside the ``conormal'' vector $\kappa\eta$. The term $b_j(\vec{y})$ describes the tracer coupling on the boundary. 

As to boundary conditions $b_1,b_2$ for the two equations of the phosphorus cycle, \cite{pa05} indicate that ``any remaining particulate organic matter that reaches the bottom of the model domain is instantly re\-mineralized''. Additionally, the model equations contain no sources or sinks for phosphate and $DOP$. This means that the total concentration (or mass) of the first two tracers remains constant with respect to time. We formalize the time-dependent total mass by 
 \begin{equation}\label{eq:defmass}
\text{mass}(y_1,y_2):=\int_{\Omega}(y_1+y_2)dx.
\end{equation}
Mathematically, conservation of mass thus corresponds to the condition
 \begin{displaymath}
\frac{\text{d}}{\text{dt}}\text{mass}(y_1,y_2)=\int_{\Omega}\partial_t(y_1+y_2)dx=0\quad\text{for all $t$}.
\end{displaymath}
On this basis we are able to derive boundary conditions. With the help of the model equations we replace the sum of the partial derivatives with respect to time $\partial_t(y_1+y_2)$ and obtain
\begin{align}\label{eq:bound1}
0=&\int_\Omega \textrm{div}(\vec{v}(y_1+y_2)-\kappa\nabla (y_1+y_2))dx\\
&+\int_{\Omega_1}(1-\nu)G(y_1,y_3)dx-\int_{\Omega_2}(1-\nu)\int_0^{h_e}G(y_1,y_3)dx_3\frac{b}{\bar{h}_e}\left(\frac{x_3}{\bar{h}_e}\right)^{-b-1}dx\nonumber.
\end{align}
Referring to the hindmost term as $M$, we obtain by inserting the definition of $\Omega_2$ and solving the emerging integral with respect to $x_3$ analytically
\begin{align*}
M:
&=(1-\nu)\int_{\Omega^{\prime}_2}\!\int_0^{h_e}\!\!G(y_1,y_3)dx_3\frac{b}{\bar{h}_e}\int_{\bar{h}_e}^{h(x^{\prime})}\!\left(\frac{x_3}{\bar{h}_e}\right)^{-b-1}\!\!dx_3dx^{\prime\!\!}\\
&=(1-\nu)\int_{\Omega^{\prime}_2}\!\int_0^{h_e}\!\!G(y_1,y_3)dx_3
\left(\frac{1}{\bar{h}_e}\right)^{-b}\left[-x_3^{-b}\right]_{\bar{h}_e}^{h(x^{\prime})}dx^{\prime\!\!}\\
&=(1-\nu)\int_{\Omega^{\prime}_2}\!\int_0^{h_e}\!\!G(y_1,y_3)dx_3
\left(1-\left(\frac{h(x^{\prime})}{\bar{h}_e}\right)^{-b}\right)
dx^{\prime\!\!}.
\end{align*}
Taking into account the definition of $\Omega_1
$, 
the integral with respect to $\Omega_1$ in Eq.~\eqref{eq:bound1} equals 
\begin{align*}
\int_{\Omega_1}(1-\nu)G(y_1,y_3)dx=\int_{\Omega^{\prime}}\int_0^{h_e}(1-\nu)G(y_1,y_3)dx_3dx^{\prime}.
\end{align*}
Finally, Gau\ss' divergence theorem yields for the first integral in Eq.~\eqref{eq:bound1}
\begin{align*}
\int_\Omega& \textrm{div}(\vec{v}(y_1+y_2)-\kappa\nabla (y_1+y_2))dx\\
&=\int_\Gamma (\vec{v}\cdot\eta(y_1+y_2)-\nabla (y_1+y_2)\cdot(\kappa\eta))ds
=-\int_\Gamma\nabla (y_1+y_2)\cdot(\kappa\eta)ds.
\end{align*}
The last equality sign holds because of the assumptions about $\vec{v}$. Combining the results obtained and heeding $\Omega^{\prime}=\Omega^{\prime}_1\,\dot\cup\,\Omega^{\prime}_2$, we transform Eq.~\eqref{eq:bound1} into
\begin{align*}
\int_\Gamma&\nabla (y_1+y_2)\cdot(\kappa\eta)ds\\
&=(1-\nu)\left(\int_{\Omega^{\prime}_1}\int_0^{h_e}G(y_1,y_3)dx_3dx^{\prime}+\int_{\Omega^{\prime}_2}\!\int_0^{h_e}\!\!G(y_1,y_3)dx_3
\left(\frac{h(x^{\prime})}{\bar{h}_e}\right)^{-b}
dx^{\prime\!}\right).
\end{align*}
Since the integral over $\Omega_j^{\prime}$ corresponds to the boundary integral over $\Gamma_j$
, the last result allows us to formulate proper boundary conditions. First of all, the right-hand side is independent of $y_2$ such that a reasonable definition is 
$b_2(\vec{y})=0 \text{ on }\Gamma.$ In accordance with the integrands of the above integral condition, we define the boundary coupling term for the first equation by
\begin{displaymath}
b_1(\vec{y}):=
\begin{cases}
-(1-\nu)\int_0^{h_e}G(y_1,y_3)dx_3&\mbox{in $\Gamma_1$,}\\
-(1-\nu)\int_0^{h_e}G(y_1,y_3)dx_3\left(\frac{h(x^{\prime})}{\bar{h}_e}\right)^{-b}&\mbox{in $\Gamma_2$}\\
\end{cases}
\end{displaymath}
as well as $b_1(\vec{y})=0$ on $\Gamma^{\prime}$.

The third tracer's behavior on the boundary remains unspecified by \cite{pa05}. Since there is a source of iron, it is not appropriate to claim conservation of mass. For this reason, the boundary conditions are limited to expressing that iron does not escape through the boundary. Speaking in mathematical terms, we impose the homogeneous boundary condition $b_3(\vec{y})=0$ on $\Gamma$.

\subsection{A weak formulation}\label{sec:weak}
The analysis in this study refers to weak solutions. These solve a ``weak formulation'', derived from the classical model equations (cf. Eq.~\eqref{eq:zustand}) by relaxing the regularity requirements.
%
%
%
Weak formulations are often preferred to classical ones because they are defined on Hilbert spaces and thus a well investigated theory about existence and uniqueness of solutions is available. The consideration of weak solutions is justified by the fact that every weak solution automatically solves the classical formulation, provided that the corresponding regularity assumptions are fulfilled. 

An important area of application is optimal control theory which is designed to solve optimality problems governed by weak formulations of differential equations.
Since parameter identification leads to such a problem, weak solutions will be required in the associated Sec.~\ref{sec:paropt}.

A weak formulation is derived by integrating the original equation multiplied with a ``test function''. By integration by parts, regularity is partly transferred to the test function. An appropriate weak formulation for the $N$-$DOP$-$Fe$ model equations is
\begin{align*}\label{eq:schwach_speziell}
&\int_0^T\{\langle y_1^{\prime},w_1\rangle +B(y_1,w_1) +\int_\Omega(-\lambda y_2+d_1(\vec{y}))w_1dx
 + \int_\Gamma b_1(\vec{y})w_1ds\}dt=0\\
 &\int_0^T\{\langle y_2^{\prime},w_2\rangle +B(y_2,w_2) +\int_\Omega(\lambda y_2+d_2(\vec{y}))w_2dx
 \}dt=0\\
 &\int_0^T\{\langle y_3^{\prime},w_3\rangle +B(y_3,w_3) +\int_\Omega(J_{Fe}(y_3)-\lambda y_2R_{Fe}+ d_3(\vec{y}))w_3dx
 \}dt\\
 &\qquad\qquad\qquad\qquad\qquad\qquad\qquad\qquad\qquad\qquad\qquad\qquad=\int_0^T\int_\Omega S_{Fe}w_3dxdt.
 \end{align*}
By $w_1,w_2,w_3$ we denote the test functions. The angle brackets indicate that the temporal derivatives $y_j^{\prime}$ are elements of a ``dual space'', i.e. maps on the space of test functions. This kind of differentiability, usually referred to as ``distributional'' \citep[Def.~23.15]{zei90}, is weaker than differentiability in the usual sense. The time-dependent term $B$, defined by 
\begin{equation*}
 B(y_j,w_j):=\int_\Omega(\kappa\nabla y_j\cdot\nabla w_j)dx+  \int_\Omega\textrm{div}(\vec{v} y_j)w_jdx
 \end{equation*}
 summarizes the ocean transport. The first summand is the result of applying integration by parts or, more precisely, Green's identity \citep[Chap.~7, Eq.~(G1)]{str08} to the diffusion term. In this connection, the boundary conditions, formerly formulated on their own, enter the weak formulation.  

Leaving out test functions and integrals, the weak formulation becomes the equivalent ``operator equation''
\begin{align}\label{eq:problemdualraum}
 y_1^{\prime}+B(y_1)-\lambda y_2&=-d_1(\vec{y})-b_1(\vec{y})\nonumber\\
  y_2^{\prime}+B(y_2)+\lambda y_2&=-d_2(\vec{y})\\
  y_3^{\prime}+B(y_3)+J_{Fe}(y_3)-\lambda y_2R_{Fe}&=S_{Fe}-d_3(\vec{y})\nonumber,
\end{align}
valid in the already mentioned dual space. In the following, we deal with the adjusted version of Eq.~\eqref{eq:problemdualraum} distinguished by the appearance of $J$ instead of $J_{Fe}$ in the last equation.

The weak formulation introduced here is derived in detail by \cite{ro14}. Further information about weak formulations and operator equations can be found e.g. in \citet{tr02,gaj67,zei90}.

Three different types of solution (transient, periodic, stationary) will be dealt with successively in the following sections.  

\section{Transient solutions}
This section deals with solutions with a given initial concentration. Beside the weak formulation, they satisfy the initial value condition
$y_j(x,0)=y^j_0(x) \text{ for $x\in \Omega$}$
with a prescribed concentration $y_0^j$ at the point of time $t=0$. Solutions of this kind are called transient.

The adjusted $N$-$DOP$-$Fe$ model equations have a unique transient solution for every at least quadratically integrable initial concentration. Furthermore, the solution depends continuously on the initial concentration.

To justify this assertion we refer to our recent analysis of systems analogous to Eq.~\eqref{eq:problemdualraum} with an arbitrary number of equations \citep{ro14}. The main challenge is the treatment of the nonlinear coupling terms. In our analysis, we followed \citet[Sect.~9.2, Thm.~2]{evans} by converting the problem of solving the weak formulation into a fixed point problem. To achieve this, an arbitrary fixed vector $\vec{z}=(z_1,z_2,z_3)$ is inserted into the reaction terms leading to the problem 
\begin{align*}
 y_1^{\prime}+B(y_1)&=\lambda z_2-d_1(\vec{z})-b_1(\vec{z})\\
  y_2^{\prime}+B(y_2)&=-\lambda z_2-d_2(\vec{z})\\
  y_3^{\prime}+B(y_3)+J(y_3)&=S_{Fe}+\lambda z_2R_{Fe}-d_3(\vec{z})\\
 y_j(0)&=y^j_0 \text{ for all $j\in \{1,2,3\}$.}
\end{align*}
Unique solvability of linear or monotone equations is a well-investigated problem \citep{zei90,la68,tr02}. While the first two equations are linear, the third is monotone because the reaction term $J$ is based on a monotonically increasing real function (cf. Sec.~\ref{subsec:iron}). In addition, both $B$ and $J$ have a linear growth and $B$ is monotone and linear. Therefore, the system depending on $\vec{z}$ has a unique solution $y$. 

Obviously, 
a fixed point of the map $\vec{z}\mapsto \vec{y}=(y_1,y_2,y_3)$ corresponds to a transient solution. Banach's Fixed Point Theorem \citep[Thm.~1.A]{zei85} yields the existence of a unique fixed point provided that the reaction terms on the right-hand side are Lipschitz continuous with a Lip\-schitz constant independent of time. In the $N$-$DOP$-$Fe$ model, these terms are either linear (multiplication with $\lambda$) or based on the uptake function $G$. However, $G$ proves Lipschitz continuous on $\mathbb{R}^2$, equipped with an arbitrary norm $\|.\|$, according to the estimation
\begin{align*}
|G(y_1,&y_3)-G(w_1,w_3)|\leq
\alpha\left|\frac{y_1}{|y_1|+K_P}\frac{y_3}{|y_3|+K_F}
-
\frac{w_1}{|w_1|+K_P}\frac{w_3}{|w_3|+K_F}\right|\\
&\leq
\alpha\left(\left|\frac{y_1}{|y_1|+K_P}-\frac{w_1}{|w_1|+K_P}\right|\left|\frac{y_3}{|y_3|+K_F}\right|\right.\\
&\qquad\qquad\qquad\qquad\qquad\qquad\quad+
\left.\left|\frac{w_1}{|w_1|+K_P}\right|\left|\frac{y_3}{|y_3|+K_F}-\frac{w_3}{|w_3|+K_F}\right|\right)\\
&\leq
\alpha\left(\frac{1}{K_P}\left|y_1-w_1\right|
+
\frac{1}{K_F}\left|y_3-w_3\right|\right)
\leq
C \|(y_1,y_3)-(w_1,w_3)\|,
\end{align*}
using that saturation functions are Lipschitz continuous and that the Lipschitz constant equals the inverse of the half saturation constant. Furthermore, the absolute value of the saturation functions is bounded by one. $C$ is a constant depending on $\alpha,K_P, K_F$.
%
%
%
%
%

It can be stated as a general observation that model equations have transient solutions if all reaction terms are either Lipschitz continuous with a Lipschitz constant independent of time or if they are monotone with linear growth.
\section{Periodic solutions}\label{sec:periodic}
In this section, we consider the existence of steady annual cycles or periodic solutions. These are characterized 
by the additional condition
 $y_j(x,0)=y_j(x,T) \text{ for $x\in \Omega$}$
signifying that the solution reaches its initial value again at the time $T$. Interpreting $T$ as one year,
periodic solutions correspond to steady annual cycles, generally approximated via spin-up or fixed-point iteration. 
Unfortunately, it is seldom possible to prove analytically, that the fixed-point iteration converges. 

Given a fixed $C>0$, there is a periodic solution of the adjusted $N$-$DOP$-$Fe$ model with the additional property that the sum of the first two tracers' mass (hereafter called total mass although iron is excluded) is equal to $C$. In particular, a periodic solution is not unique. The time-dependent total mass is formalized in Eq.~\eqref{eq:defmass}. Remark that, due to the choice of reaction terms and boundary conditions, every solution has a constant total mass with respect to time.

The justification of the existence result bases on our previous study about periodic solvability of the $N$-$DOP$ model, assuming a constant distribution of iron \citep{ro214}.  We will present our proceeding and indicate how to extend the argumentation to the model with three equations.

In opposite to the transient case, the model equations are simplified by inserting a fixed $\vec{z}$ only into the reaction terms $d_j$ and $b_j$. In a first step, we solve the periodic problem
\begin{align*}
 y_1^{\prime}+B(y_1)-\lambda y_2&=-d_1(\vec{z})-b_1(\vec{z})\\
  y_2^{\prime}+B(y_2)+\lambda y_2&=-d_2(\vec{z})\\
  y_3^{\prime}+B(y_3)+J(y_3)-\lambda y_2R_{Fe}&=S_{Fe}-d_3(\vec{z})\\
 y_j(0)&=y_j(T) \text{ for all $j\in \{1,2,3\}$}
\end{align*}
and find a fixed point of the map $\vec{z}\mapsto \vec{y}$ afterwards. Since there is no standard technique for periodic solvability of such systems, we had to develop a method of our own. Its presentation provides insight into the behavior of closed systems.

Being partly decoupled, the three equations of the system above can be solved consecutively. The second equation has a unique periodic solution depending on $\vec{z}$ because 
the linear summand $+\lambda y_2$ allows the estimate
\begin{equation}\label{eq:coercive}
\int_0^T\int_\Omega \lambda y_2y_2dxdt\geq \lambda\|y_2\|^2
\end{equation}
where $\|.\|$ is the norm in a Hilbert space of functions depending on space and time. 
%
%
This condition causes the corresponding fixed-point iteration to converge. 
 Since the first equation lacks a comparable summand, we first solve $S^{\prime}+B(S)=-d_1(\vec{z})-b_1(\vec{z})-d_2(\vec{z})$, the equation for the sum $S=y_1+y_2$. A standard existence theorem provides a unique solution in a special function space with the additional condition $\int_\Omega Sdx=0$ \citep[Chap.~VI, Thm.~1.4]{gaj67}. 
The properly chosen boundary conditions (and therefore conservation of mass) ensure that $S$ belongs to the desired function space. Due to the choice of $S$, the function $y_1=S-y_2+|\Omega|^{-1}C$ solves the first equation. Here, the symbol $|\Omega|$ stands for the measure of the domain $\Omega$. The solution $\vec{y}=(y_1,y_2)$ is unique and satisfies $\text{mass}(y_1,y_2)=\int_\Omega|\Omega|^{-1}Cdx=C$.
%
%
%

%
The third equation can be solved using the mentioned standard theorem as well. 
The adjusted coupling term $J$ is strictly monotone because $FeF$ has a strictly increasing slope (cf. Sec.~\ref{subsec:iron}) and all other factors are positive. Additionally assuming a positive lower threshold $c_p>0$ for the particle concentration $C_p$, 
the term $J$ satisfies
\begin{equation}
\label{eq:coercive2}
\int_0^T \int_\Omega J(y_3)y_3dxdt\geq C_1\|y_3\|^2-C_2\|y_3\|.
\end{equation} 
This condition is slightly weaker than the one displayed in Eq.~\eqref{eq:coercive} but still sufficient to ensure that the operator $B+J$ is ``coercive'' which is a crucial property in the existence theorem. Actually, the original reaction term $J_{Fe}$ was replaced because it lacks coercivity. Taking into account the properties of $J$, the standard theorem yields a unique periodic solution $y_3$ depending on $\vec{z}$.

As in the transient case, 
a fixed point of the map $\vec{z}\mapsto \vec{y}=(y_1,y_2,y_3)$ corresponds to a periodic solution of the original problem. Schauder's Fixed Point Theorem \citep[Thm.~2.A]{zei85} yields the desired fixed point because the coupling terms $d_j$ and $b_j$ are bounded independently of $\vec{z}$. This is due to the uptake $G(z_1,z_3)$ whose absolute value is bounded by 1. 

As a general observation, we can state that a proof of periodic solvability requires a specialized proceeding adapted to the model in question. Concerning the choice of reaction terms, it can be said that every function bounded independently of the inserted argument is allowed. Furthermore, the $j$-th equation can (and to a certain extent has to) contain reaction terms which are monotone and coercive, i.e. satisfy an estimate analogous to Eq.~\eqref{eq:coercive2}, and depend only on the $j$-th tracer.

\section{Stationary solutions}
Stationary solutions represent the constant tracer concentrations reached with a fixed forcing. Consequently, they solve a time-independent variant of the original model Eq.~\eqref{eq:zustand}, characterized by vanishing temporal derivatives and temporally constant terms $\vec{v},\kappa,d_j,b_j,J_{Fe}$. The corresponding adjusted weak formulation including test functions is
\begin{align}\label{eq:problemstationary}
&B(y_1,w_1) +\int_\Omega(-\lambda y_2+d_1(\vec{y}))w_1dx
 + \int_\Gamma b_1(\vec{y})w_1ds=0\nonumber\\
 &B(y_2,w_2) +\int_\Omega(\lambda y_2+d_2(\vec{y}))w_2dx
 =0\\
 &B(y_3,w_3) +\int_\Omega(J(y_3)-\lambda y_2R_{Fe}+ d_3(\vec{y}))w_3dx
 =\int_\Omega S_{Fe}w_3dx.\nonumber
 \end{align}
 Since the equation is independent of time, the integrals over $[0,T]$ are missing. For the same reason, the test functions, like the solutions, only depend on the spatial coordinates.

It is important to bear in mind that stationary solutions are not periodic in the sense of the last section although they are constant with respect to time and thus initial and terminal values coincide. The reason is that stationary solutions in the sense of this study correspond to a constant forcing, while periodic solutions solve the equations with temporally variable summands. Thus, the periodic solutions of the last section are not necessarily constant.

The adjusted $N$-$DOP$-$Fe$ model equations have a stationary solution, i.e. a solution of Eq.~\eqref{eq:problemstationary}, for every total mass $C>0$. 

An analysis of stationary solutions has not been published yet. However, since they are closely related to periodic solutions, the considerations of the last section can be transferred. 
As in the periodic case, the system in Eq.~\eqref{eq:problemstationary} is simplified by inserting a fixed $\vec{z}$ into the reaction terms $d_j$ and $b_j$. The equations are solved in the same order as the periodic ones, now with the help of the theorem of Browder and Minty \citep[Thm.~26.A]{zei89} instead of the periodic standard theorem.
The theorems are closely related. Both claim the assumptions continuity, monotonicity and coercivity all of which have already been investigated in the last section. Finally, Schauder's Fixed Point Theorem is applied in exactly the same way as before.

Due to the analogous proceeding, the conclusions concerning the choice of reaction terms are the same as in the last section.

\section{Parameters in the $N$-$DOP$ model}\label{sec:paropt}
%
Parameter identification is one of the most important areas of application of the $N$-$DOP$ model. 
As explained in the introduction, parameters are identified via minimizing the difference between model output, depending on the parameters, and original data. Numerical optimization methods to be tested are applied to synthetic data and judged by their ability to identify the corresponding known optimal parameters. 

A criterion for the explanatory power of this approach is the parameters' unique identifiability. Parameters are called uniquely identifiable if each set of parameter values is associated with a single model output. Otherwise, the parameters are called dependent.

Thus, information about dependencies in the $N$-$DOP$ model contributes to the validation of numerical tests and the interpretation of results like those of \cite{pr13} whose method did not identify all parameters correctly.

%
%
\subsection{Investigation of identifiable and dependent parameters
}\label{subsec:dep}

In this section, we investigate the seven $N$-$DOP$ model parameters as to possible dependencies. As a result, we will be able to tell which parameters are uniquely identifiable and thus suited for testing purposes and which are (supposedly) dependent. 

To this end, we assume that the equations associated with the parameter sets $\vec{u}_1=(\lambda_1,\alpha_1,K_{P1},K_{I1},K_{W1},b_1,\nu_1)$ and $\vec{u}_2=(\lambda_2,\alpha_2,K_{P2},K_{I2},K_{W2},b_2,\nu_2)$ have the same nontrivial solution $\vec{y}=(y_1,y_2)$. It is not relevant if the solution is periodic, transient or stationary. Parameters that prove equal in both parameter vectors are uniquely identifiable. We limit the investigation to the natural case that $\vec{y}$ has two nontrivial components because otherwise there are obvious dependencies. For instance, the solution $y_1=0$ allows arbitrary parameters $\alpha,K_{P},K_{I},K_{W}$ and $y_2=0$ allows an arbitrary $\lambda$. %
%
%
%

It is reasonable to restrict the considerations to positive $\alpha$, $\nu<1$ and $K_I>0$. 
The value $\alpha=0$ is not a likely maximum production and obviously effects dependencies since $\alpha=0$ and $K_I,K_P,K_W$ arbitrarily chosen lead to the same solution. In case $\nu=1$, the export into the aphotic zone is zero because all consumed phosphate is transformed into $DOP$ and thus the sinking parameter $b$ can be chosen arbitrarily. Finally, a vanishing half saturation constant $K_I=0$ eliminates the influence of light. All other parameters are assumed to be nonnegative.  

The choice of $\vec{y}$ signifies that
\begin{align*}
 y_1^{\prime}+B(y_1)&=\lambda_i y_2-d_1(\vec{u}_i,\vec{y})-b_1(\vec{u}_i,\vec{y})\\
  y_2^{\prime}+B(y_2)&=-\lambda_i y_2-d_2(\vec{u}_i,\vec{y})
  \end{align*}
hold for both $i=1,2$. This time, we explicitly indicate the reaction terms' dependence on the parameter vectors. Clearly, the left-hand sides are equal 
for both $i$ and therefore also the right-hand sides, i.e.
\begin{align}
\lambda_1 y_2-d_1(\vec{u}_1,\vec{y})-b_1(\vec{u}_1,\vec{y})&=\lambda_2 y_2-d_1(\vec{u}_2,\vec{y})-b_1(\vec{u}_2,\vec{y})\label{eq:pargle1}\\
-\lambda_1y_2-d_2(\vec{u}_1,\vec{y})&=-\lambda_2 y_2-d_2(\vec{u}_2,\vec{y}).\label{eq:pargle2}
\end{align}
In order to draw conclusions about the parameters, we utilize the reaction terms' specifications in each part of $\Omega$ (see Sec.~\ref{sec:phosphor}). Equation~\eqref{eq:pargle2}, considered in $\Omega_2$, reveals that $\lambda_1 y_2=\lambda_2 y_2$. Since $y_2$ is not trivial, this shows $\lambda_1=\lambda_2$, i.e. unique identifiability of $\lambda$. 
%

Taking into account $\lambda_1=\lambda_2$, Eq.~\eqref{eq:pargle1} in $\Omega_1$ yields 
\begin{equation}
\label{eq:gleichheitG}
G(\vec{u}_1,y_1)=G (\vec{u}_2,y_1).
\end{equation}
%

From Eq.~\eqref{eq:pargle2} in $\Omega_1$, we conclude with the help of Eq.~\eqref{eq:gleichheitG} 
\begin{align*}
0=\nu_1G(\vec{u}_1,y_1)-\nu_2G (\vec{u}_2,y_1)&=(\nu_1-\nu_2)G(\vec{u}_1,y_1) +\nu_2(G (\vec{u}_1,y_1)-G(\vec{u}_2,y_1))\\
&=(\nu_1-\nu_2)G(\vec{u}_1,y_1).
\end{align*}
Since $\alpha>0$ and $y_1$ is not trivial by assumption we obtain $G(\vec{u}_1,y_1)\neq0$ and thus $\nu_1=\nu_2$, i.e. unique identifiability of $\nu$. 

In the next step, we deal with the parameter $b$. Since otherwise this parameter would not occur in the model equations, we assume that the integral over $G(\vec{u}_1,y_1)$ with respect to depth is not zero everywhere. Taking into account the aphotic boundary $\Gamma_2$, it follows from Eq.~\eqref{eq:pargle1} 
 \begin{align*}
0&=\int_0^{h_e}G(\vec{u}_1,y_1)dx_3\left(\frac{h(x^{\prime})}{\bar{h}_e}\right)^{-b_1}
-\int_0^{h_e}G(\vec{u}_2,y_1)dx_3\left(\frac{h(x^{\prime})}{\bar{h}_e}\right)^{-b_2}\\
&=\int_0^{h_e}G(\vec{u}_1,y_1)dx_3\left(\left(\frac{h(x^{\prime})}{\bar{h}_e}\right)^{-b_1}
-\left(\frac{h(x^{\prime})}{\bar{h}_e}\right)^{-b_2}\right)\\
&\qquad\qquad\qquad\qquad\qquad+\left(\frac{h(x^{\prime})}{\bar{h}_e}\right)^{-b_2}\int_0^{h_e}(G(\vec{u}_1,y_1)-G(\vec{u}_2,y_1))dx_3
\end{align*}
due to $1-\nu\neq0$.
The last summand vanishes according to Eq.~\eqref{eq:gleichheitG}. Because of the assumption concerning the integral over $G(\vec{u}_1,y_1)$, we conclude 
\begin{equation*}
\left(\frac{h(x^{\prime})}{\bar{h}_e}\right)^{-b_1}
-\left(\frac{h(x^{\prime})}{\bar{h}_e}\right)^{-b_2}=0.
\end{equation*}
The fraction $q:=h(x^{\prime})/\bar{h}_e$ is strictly less than 1 for at least one $x^{\prime}$ since otherwise the domain would lack the aphotic part and the corresponding boundary. Since the natural logarithm $\ln:\mathbb{R}_{>0}\to\mathbb{R}$ is bijective we obtain 
  \begin{equation*}
-b_1\ln(q)
=\ln(q^{-b_1})
=\ln(q^{-b_2})=-b_2\ln(q).
\end{equation*} 
As ensured before, $\ln(q)\neq0$. It follows $b_1=b_2$.

The remaining parameters are incorporated in the uptake function $G$. Equation~\eqref{eq:gleichheitG} states that 
\begin{align*}
\alpha_1\frac{y_1}{|y_1|+K_{P1}}\frac{I\text{e}^{-x_3K_{W1}}}{|I\text{e}^{-x_3K_{W1}}|+K_{I1}}-
 \alpha_2\frac{y_1}{|y_1|+K_{P2}}\frac{I\text{e}^{-x_3K_{W2}}}{|I\text{e}^{-x_3K_{W2}}|+K_{I2}}=0.
\end{align*}
This equality provides the possibility to derive a condition for dependencies. The computation provided in App.~\ref{app} shows that dependencies exist if and only if it is possible to depict the solution of the first equation by
\begin{align*}
|y_1|=
\frac{c_1I+c_2\text{e}^{x_3c_7}
-c_3\text{e}^{x_3c_8}}
{c_4I-c_5\text{e}^{x_3c_7}
+c_6\text{e}^{x_3c_8}}
\end{align*}
with constants $c_1,\dotsc,c_8$ satisfying the additional condition  
$c_1={c_2}/{c_5}-(c_4+1){c_3}/{c_6}$.
In other words, if the absolute value of $y_1$ can be expressed in the indicated manner, it is possible to find two parameter vectors $\vec{u}_1,\vec{u}_2$ leading to the same solution. We cannot rule out that this condition holds for some coefficients $c_1,\dotsc,c_8$. Thus, we have to assume that dependencies exist and the four parameters $\alpha,K_{P},K_{I},K_{W}$ may not be uniquely identifiable.

This result coincides with the finding of \cite{pr13} who were able to identify $\lambda,b$ and $\nu$ (which they call $\sigma$) very well. The results for the four other parameters were less satisfying. Whereas the optimal $\alpha$ was almost reached, the approximation hardly moved towards the optimal $K_{I}$. Concerning the other parameters, the accuracy of approximation lay in between.

\subsection{Elimination of dependencies}
To be valid and significant, tests of numerical methods should rely on models containing only identifiable parameters. In this section, we therefore propose two possible ways to eliminate the $N$-$DOP$ model's dependencies.

The first possibility is to identify a reduced number of parameters after fixing the remaining ones. The fixed values have to be quantified in another way (experiments, estimations). 
In case of the $N$-$DOP$ model we propose to fix $K_I$ and $K_W$. Then the five parameters $\lambda,\alpha,K_{P},b,\nu$ remain to be identified via optimization. 

We justify this suggestion by proving the unique identifiability of the five variable parameters. This property has already been shown for $\lambda,b,\nu$. Concerning the two remaining parameters, the analog of Eq.~\eqref{eq:gleichheitG},
\begin{align*}
\alpha_1\frac{y_1}{|y_1|+K_{P1}}=
 \alpha_2\frac{y_1}{|y_1|+K_{P2}},
\end{align*}
can be simply transformed into
\begin{align*}
\frac{\alpha_1}{\alpha_2}=
 \frac{|y_1|+K_{P1}}{|y_1|+K_{P2}}=1+\frac{K_{P1}-K_{P2}}{|y_1|+K_{P2}}.
\end{align*}
The left-hand side is clearly constant. Since the only constant solution $y_1=0$ is excluded, the right-hand side is constant if and only if the last fraction vanishes, i.e. if $K_{P1}=K_{P2}$. We directly conclude $\alpha_1/\alpha_2=1$, i.e. $\alpha_1=\alpha_2$. Thus, both remaining parameters are identifiable.

We state 
that 
slightly reducing the number of unknown parameters in the $N$-$DOP$ model eliminates all dependencies. 
Thus, the $N$-$DOP$ model with the reduced parameter set is recommendable to test and validate numerical methods reliably.

A second possibility to construct a model with identifiable parameters is to replace the reaction terms being responsible for dependencies. This may involve a reduced number of parameters which are not immediately interpretable with respect to the biochemical processes modeled.

\begin{figure}
\includegraphics[width=8.3cm]{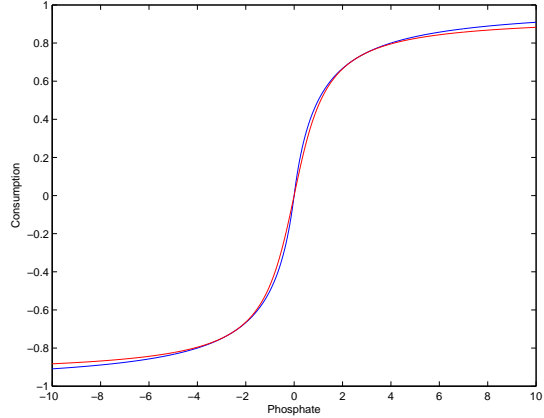}
\caption{Comparison of saturation function (blue) with $\alpha=K=1$ and the arc tangent scaled with the parameter $\beta=0.6$ (red).}
\label{fig:arctan}
\end{figure}

In case of the $N$-$DOP$ model, the results indicate that dependencies are caused by the product of saturation functions. A possible alternative for the saturation function $\alpha y_1/(|y_1|+K)$ with two parameters $\alpha,K$ could be the scaled arc tangent $\beta\arctan(y_1)$ with only one parameter $\beta$. 
As we see in Fig.~\ref{fig:arctan}, the curves of both functions with appropriate parameter values behave similarly.

Alternatives to the arc tangent are for example an approximation via the Taylor series or a linear function. Which one is suited best, depends e.g. on the values the solution $y_1$ reaches. 
%
%
%
%
\section{Discussion}
A numerical treatment of partial differential equations is usually preceded by a mathematical ana\-lysis. Numerical methods are designed to approximate exact solutions. Thus, if there is no exact solution, the numerically obtained results are dubious and unpredictable. A mathematical analysis therefore contributes to the validation and assessment of ecosystem models.
In this work, we provide correspondent information about the $N$-$DOP$ model by \cite{pa05} which is extensively used. 

In the first part, the model formulation is stated in full mathematical detail.   
We additionally develop boundary conditions assuming conservation of mass. This condition is both common regarding ecosystem models and crucial in the proofs of periodic and stationary solvability.
Furthermore, we explicitly derive the reaction term modeling complexation and scavenging of iron 
%
%
and find an alternative formulation with a similar behavior and convenient mathematical properties. 
Finally, a weak formulation of the model equations is specified.
The exact mathematical formulation developed here is an essential premise for the following analysis.

As a result, the analysis yielded that the adjusted $N$-$DOP$-$Fe$ model and therefore, in particular, the (original) $N$-$DOP$ model has transient, stationary and periodic solutions.  

Periodic solutions or steady annual cycles, characterized by equal initial and terminal values (e.g. at the beginning and the end of one year) are of particular interest as they are required in most applications. The analysis of periodic solutions provides an interesting insight into systems conserving mass. Assuming conservation of mass, which is commonly done in the context of ecosystem models, we were able to show that there exists a distinct periodic solution for each mass contained in the phosphorus cycle. 
In particular, there are other periodic solutions than the trivial $\vec{y}=(0,0,0)$, corresponding to the mass zero. 
This result confirms the observation made in spin-up computations that every initial mass leads to a corresponding periodic solution.
%
%
%

 As a conclusion, the results suggest that it is possible and meaningful to solve the  $N$-$DOP$ model equations  numerically. 
 
Furthermore, the mathematical conditions found during the analysis enables the identification and elimination of problematic reaction terms. The reaction term for scavenging and complexation of iron $J_{Fe}$, for instance, could be substituted by a suitable alternative $J$. 

%
In the second part, we investigated the $N$-$DOP$ model's parameters. 
Four of the seven parameters proved probably dependent and thus not uniquely identifiable. This result suggests that unsatisfying results in parameter identification might be due to the model itself instead of an inadequate numerical method. As a consequence, the original $N$-$DOP$ model seems unsuited for tests of parameter identification methods.

The analysis reveals that dependencies originate from the modeling of the biological uptake limited by iron, phosphate and light. Thus, a model without dependencies either lacks this term altogether or has a reduced number of parameters to be identified. We showed that one possibility is to fix the two parameters $K_I$ and $K_W$ modeling insolation and to identify only the remaining five which are then uniquely identifiable. Especially in the context of numerical tests, this is a convenient procedure since the model itself remains unaltered.
Alternatively, by replacing the critical reaction terms altogether, new parameters without dependencies can be introduced. The analysis allowed us to formulate different substitutes in such a way that the reaction terms' mathematical behavior remains similar to the original one. Other alternatives could be developed.  As a future task, it remains to find out which alternative is in line with the ecosystem modeled. 
%

As a conclusion, the analysis improves the interpretation of tests in the context of parameter identification and the assessment of methods tested.
By using one of our proposed alternatives instead of the original $N$-$DOP$ model, one important source of uncertainty can be eliminated.
 
Finally, the analysis concerning identifiability and solvability of the $N$-$DOP$-$Fe$ model is performed with the help of universal mathematical methods. 
 Thus, the analysis conducted for the $N$-$DOP$-$Fe$ model provides the basis for the assessment of other models as well.

\section*{Acknowledgements}
The research of Christina Roschat was supported by the DFG Cluster of Excellence Future Ocean, Grant No. CP1338.

\appendix
\section{Derivation of a characterization for dependencies}\label{app}    
In order to find a characterization for dependencies, we solve the equation
\begin{align*}
\alpha_1\frac{y_1}{|y_1|+K_{P1}}\{I\text{e}^{-x_3K_{W1}}{|I\text{e}^{-x_3K_{W1}}|+K_{I1}}=
 \alpha_2\frac{y_1}{|y_1|+K_{P2}}\frac{I\text{e}^{-x_3K_{W2}}}{|I\text{e}^{-x_3K_{W2}}|+K_{I2}}
\end{align*}
for $|y_1|$. Taking into account $\alpha_i,I\text{e}^{-x_3K_{W1}}>0$ and $y_1\neq0$, we obtain
 \begin{align*}
\frac{\alpha_1}{\alpha_2}&=\frac{|y_1|+K_{P1}}{y_1}\frac{I\text{e}^{-x_3K_{W1}}+K_{I1}}{I\text{e}^{-x_3K_{W1}}}
\frac{y_1}{|y_1|+K_{P2}}\frac{I\text{e}^{-x_3K_{W2}}}{I\text{e}^{-x_3K_{W2}}+K_{I2}}\\
&=
\frac{|y_1|+K_{P1}}{|y_1|+K_{P2}}
\frac{I\text{e}^{-x_3K_{W1}}+K_{I1}}{I\text{e}^{-x_3K_{W2}}+K_{I2}}
\text{e}^{-x_3(K_{W2}-K_{W1})}.
\end{align*}
Using the abbreviation $C:=\alpha_1/\alpha_2$, we calculate
\begin{align*}
0&=
C{(|y_1|+K_{P2})}{(I\text{e}^{-x_3K_{W2}}+K_{I2})}-
{(|y_1|+K_{P1})}
{(I\text{e}^{-x_3K_{W1}}+K_{I1})}
\text{e}^{-x_3(K_{W2}-K_{W1})}\\
&=
C|y_1|(I\text{e}^{-x_3K_{W2}}+K_{I2})+CK_{P2}(I\text{e}^{-x_3K_{W2}}+K_{I2})\\
&\qquad-
|y_1|(I\text{e}^{-x_3K_{W2}}+K_{I1}\text{e}^{-x_3(K_{W2}-K_{W1})})
-K_{P1}(I\text{e}^{-x_3K_{W2}}+K_{I1}\text{e}^{-x_3(K_{W2}-K_{W1})})\\
&=
|y_1|\{C(I\text{e}^{-x_3K_{W2}}+K_{I2})
-
(I\text{e}^{-x_3K_{W2}}+K_{I1}\text{e}^{-x_3(K_{W2}-K_{W1})})\}\\
&\qquad+CK_{P2}(I\text{e}^{-x_3K_{W2}}+K_{I2})
-K_{P1}(I\text{e}^{-x_3K_{W2}}+K_{I1}\text{e}^{-x_3(K_{W2}-K_{W1})}).
\end{align*}
Rearranging the summands leads to
\begin{align*}
|y_1|
\{((C-1)I-&K_{I1}\text{e}^{x_3K_{W1}})\text{e}^{-x_3K_{W2}}
+CK_{I2}\}\\
&=
((K_{P1}-CK_{P2})I+K_{P1}K_{I1}\text{e}^{x_3K_{W1}}
)\text{e}^{-x_3K_{W2}}-CK_{P2}K_{I2}.
\end{align*}
Taking into account $K_I\neq0$, it is possible to prove that 
dependencies exist if and only if the expression in curly brackets is not zero. 
Dividing the equation by this expression and canceling out $\text{e}^{-x_3K_{W2}}$ on the right-hand side, we obtain 
\begin{equation}\label{eq:erg}
|y_1|
=\frac{(K_{P1}-CK_{P2})I+K_{P1}K_{I1}\text{e}^{x_3K_{W1}}
-CK_{P2}K_{I2}\text{e}^{x_3K_{W2}}}
{(C-1)I-K_{I1}\text{e}^{x_3K_{W1}}
+CK_{I2}\text{e}^{x_3K_{W2}}}
.
\end{equation}
Equation~\eqref{eq:erg} corresponds to the condition specified in Sec.~\ref{subsec:dep}. The additional condition for the constants $c_1,\dotsc,c_8$ reflects the dependence of the coefficients in Eq.~\eqref{eq:erg} on each other.


\bibliography{meinebibliothek}

\end{document}